\documentclass[12pt,a4paper]{article}
\usepackage[T2A]{fontenc}
\usepackage[cp1251]{inputenc}
\usepackage{amsmath,amssymb,amsthm,amsfonts,amscd}

 \theoremstyle{plain}
\newtheorem{theorem}{Theorem}
\newtheorem*{main*}{Main Theorem}
\newtheorem*{theorem*}{Theorem}
\newtheorem{lemma}[theorem]{Lemma}

\newtheorem{corollary}[theorem]{Corollary}

\theoremstyle{definition}
\newtheorem{definition}[theorem]{Definition}
\newtheorem{remark}[theorem]{Remark}
\newtheorem*{remark*}{Remark}
\newtheorem*{example*}{Example}

\def\Z{{\Bbb Z}}
\def\R{{\Bbb R}}

\def\C{{\Bbb C}}
\def\A{{\bf A}}

\def\i{{\bf i}}

\date{}

\begin{document}
\sloppy
\title{A remark on the Hopf invariant for spherical 4-braids}
\author{Akhmet'ev P.M.}

\maketitle

\section{Introduction}

An approach by J.Wu describes homotopy groups $\pi_{n}(S^2)$ of the standard 2-sphere as isotopy classes of spherical
$n+1$--strand Brunnian braids, for more details, see f.ex. \cite{B-M-V-W}, Sequence (1.1). This approach 
is not possible for $n=3$ in the case of $4$-strand braids.
%Тo explicit relation of Brunnian $4$-stand braids with $\pi_3(S^2)$ is known.  

The homotopy group $\pi_3(S^2)$ in an infinite cyclic group, detected by the Hopf invariant
\begin{eqnarray}\label{hopfinv}
H: \pi_3(S^2) \mapsto \Z. 
\end{eqnarray}
An element of $\pi_3(S^2)$ is represented by a mapping $l:S^3 \mapsto S^2$, which is considered up to homotopy. The Hopf invariant $H(l)$
is well-defined as the integer linking number of two oriented curves $l^{-1}(a)$, $l^{-1}(b)$, where
$a,b \in S^2$ be a pair of regular points of $l$. 
The Hopf invariant is very important for applications. 

The goal of the paper is to modify the definition of Brunnian spherical 4-component braids, and to define the Hopf
invariant as a function of isotopy classes of spherical braids, which are Brunnian in a new strong sense. 
The Hamiltonian provides an elegant
method for generating simple geometrical examples of complicated braids and
links, as is presented in \cite{B}.
%The Hopf invariant is related with the "`Borromean"' 3-components (ordinary) braid, which is a Brunnian braid. 

Let us formulate the following problems:

$\bullet$ Derive applications of higher-order winding numbers to generate Hamiltonian
motion of 4 vortex in two dimensions on the sphere. For 3 vortex on the plane this is done in \cite{B}.

$\bullet$ To unify the approach \cite{V} Ch.3 to $\pi_{\ast}(S^2)$ with the Wu's approach.
%This exactly means that the Borromean braids consists of 3 strings, each pair is unknotted. 

$\bullet$ To investigate integer lifts of the generator of $\pi_4(S^2)$ (the Arf invariant) by the Wu's approach. 

The paper is organized as following. In Section 2 we recall required definitions concerning first-order stage of the construction and determine the linking numbers of spherical 4-component braids. In Section 3 the Hopf invariant for 4-component spherical braids is well-defined. This is a second-order  particular defined invariant: 
to define this invariant we should assume that the all linking numbers (there are two) of components of a spherical braid are equal to zero. 
Main results are formulated in Theorems $\ref{hopf}$, $\ref{main}$.
In Section 4 we give proofs of the main results.

In a private letter (August 2013) prof. Viktor Ginzbugr  (about a draft of the paper): "`The subject is certainly interesting..."'. I am grateful to him for the interest.

The  results was presented at International Conference "Nonlinear Equations and Complex Analysis"
in Russia (Bashkortostan, Bannoe Lake) during the period since March 18 (arrival day) till March 22 (departure day), 2013.

\section{Linking numbers for spherical braids}

By a spherical (ordered) $n$--braid we mean a collection of embeddings of the standard circles
$$f:\bigcup_{i=1}^n S^1_i \subset S^2 \times S^1, $$
where the composition of this embedding with the standard projection $S^2 \times S^1 \to S^1$ on the second factor in the target space, restricted to an arbitrary component $S^1_i$, $i=1,\dots,n$ is the identity mapping $S^1_i \to S^1$.
The space of all ordered spherical $n$-braids up to isotopy is denoted by $Br_n$. It is well-known that $Br_n$ is
a group.

For a fixed value $t \in S^1$, a braid $f \in Br_n$ intersects the level $S^2 \times t$ by an (ordered) collection of $n$ points
$\{z_1(t), \dots z_n(t)\}$. Let assume that $n=4$.  Denote by
$$ g=g(f): S^{1}_{1} \cup S^1_{2} \cup S^1_{3} \subset S^2 \times S^1, $$
the $3$-component braid, obtained from $f$ by eliminating of the last component $S^1_4$.

Let us identify the sphere $S^2$ with the Riemann sphere, or with the complex projective line $\hat \C$.
For a braid $f$ let us consider the collection of M$\ddot{\rm{o}}$bius transformations, which transforms the points
$z_{1},z_{2},z_{3}$ into $0,1,\infty$ correspondingly:
$$F(z;t) = \frac{(z-z_{1}(t))(z_{2}(t)-z_{3}(t))}{(z - z_{3}(t))(z_{2}(t)-z_{1}(t))}.$$ 
The image $F(f)$ is a 4-strand braid with the constant components $\{z_{1}(t), z_{2}(t), z_{3}(t)\} = \{0,1,\infty\}$.
Denote this braid by 
\begin{eqnarray}\label{Mobius}
F(f) = f^{norm}.
\end{eqnarray}
The 3-strand braid $g$, constructed from $f^{norm}$ is the constant braid at the points $\{0,1,\infty\}$. 
The last component $f^{norm}_4(S^1)$ of $F(f)$ is represented by a closed path $z_4(t) \in \hat \C \setminus \{0,1,\infty\}, t \in \R^1/2\pi$.

For a given (ordered) 4-component braid $f$ let us define the linking number $Lk(f)$,
\begin{eqnarray}\label{link}
 Lk: Br_4 \to \Z. 
\end{eqnarray} 
Consider the following 1-form
\begin{eqnarray}\label{omega0}
\omega_0 = \frac{1}{2 \pi \i} \frac{dz}{z}. 
\end{eqnarray}
By definition we get
$$ d \log(z) = \frac{1}{2 \pi \i} \frac{dz}{z},$$ 
where $\log(z)$ is given by the formula:
$$ \log(z) = (2\pi \i)^{-1} \int  \frac{dz}{z},$$
assuming that $\log(1)=0$, as a multivalued complex function. 

Define $Lk(f)$ by the formula:
$$Lk(f) = \Re{\int_{0}^{2\pi} \frac{d z_4(t)}{z_4(t)}} = \int_{f^{norm}_4} \omega_0,$$
where $\Re$ is the real part of the integral.
By construction,  $Lk(f)$ is the winding number, i.e. the integer number of rotations of the path $z_4(t)$ with respect to the origin and the infinity in $\hat\C$.

The permutation group $\Sigma(4)$ of the order $24$ acts on the space of ordered spherical braids:
\begin{eqnarray}\label{action}
 \Sigma(4) \times Br_4 \to Br_4. 
 \end{eqnarray}
The image of an ordered braid
$f$ by a transposition $\sigma:(1,2,3,4)\mapsto(\sigma_1,\sigma_2,\sigma_3,\sigma_4)$ is well-defined by the corresponding re-ordering of components of $f$.
Let us investigate the orbit of the linking numbers $Lk(f)$ with respect to $(\ref{action})$. Simply say, we investigate how many independent linking numbers of components of braids are well-defined? 

Let us consider the following exact sequences of groups:
\begin{eqnarray}\label{first}
          0    \longrightarrow    \A_4  \longrightarrow    \Sigma_4 \longrightarrow \Z/2 \longrightarrow 0,
\end{eqnarray}
\begin{eqnarray}\label{second}
 0 \longrightarrow \Z/2 \times \Z/2  \longrightarrow \A_4  \longrightarrow \Z/3  \longrightarrow 0.
\end{eqnarray}
The subgroup $\A_4 \subset \Sigma_4$ in the sequence $(\ref{first})$ is represented by permutations, which preserve signs (equivalently, which is  decomposed into an even number of elementary transpositions).
The subgroup $\Z/2 \times \Z/2 \subset \A_4$ in the sequence $(\ref{second})$ is generated by the  permutations $\{(1,2)(3,4); (1,3)(2,4); (1,4)(2,3)\}$. 

Let us consider $2$-primary subgroup $K \subset \Sigma$ of the order $8$, which is defined as the extension of the subgroup $\Z/2 \times \Z/2$ from the sequence $(\ref{second})$, been included in the sequence $(\ref{first})$. An
epimorphism
\begin{eqnarray}\label{epi}
\theta = (\theta_1,\theta_2): K \to \Z/2 \times \Z/2,
\end{eqnarray} 
is defined as follows: $\theta_1(\sigma)=1$ (the group $\Z/2$ is in the multiplicative form), if $\sigma$ preserves a (non-ordered) partition $(1,3)(2,4)$, and
$\theta_1(\sigma)=-1$, otherwise. Therefore $\theta_1$ is an epimorphism with the kernel $\Z/2 \times \Z/2$ from the left subgroup of the sequence $(\ref{second})$. The epimorpism $\theta_2(\sigma)$ is determined by the sign of  a permutation $\sigma$, this is the restriction of the right epimorphism in the sequence $(\ref{first})$ to the subgroup $K \subset \Sigma_3$.
%The  epimorpism $\theta$ extends to an epimorphism $\Sigma_4 \to \Z/2 \times \Z/2$.  

\begin{lemma}\label{linking}

--1. The function $(\ref{link})$ is invariant with respect to the action $(\ref{action})$ (the re-numbering of components) by an arbitrary permutation, which  in the kernel of $\theta$ in $(\ref{epi})$,
and is skew-invariant for the action by a permutation, which is in the kernel of $\theta_1$ (the composition of $\theta$ with the projection on the first factor), but not in the kernel of $\theta_2$  (the composition of $\theta$ with the projection on the second factor).

--2. Denote by $\tilde{f} \in Br_4$ the ordered braid, which is obtained from $f \in Br_4$ by the action $(\ref{action})$  by the element $(1,2)$ ($\theta((1,2)) \in \Z/2 \times \Z/2$ is the product of the generators). There exists an ordered braids $f \in Br_4$,  for which the linking numbers $Lk(f)$, $Lk(\tilde{f})$ are arbitrary integers. 

\end{lemma}

From Lemma one may deduce the following corollary.

\begin{corollary}\label{winding}

--1. For an arbitrary braid $f \in Br_4$ the linking number $Lk(f)$ is well-defined as the differences of the winding number of the component $2$  between the components $1$ and $3$ with
the winding number of the component $4$ between the components $1$ and $3$. 

--2. For a braid $\tilde{f} \in Br_4$, where $f \in Br_4$ is an arbitrary, $\tilde{f}$ is defined in Lemma $\ref{linking}$, 
the linking number $Lk(\tilde{f})$ is well-defined as the winding number of the component $2$ between the components $1$ and $3$ with the winding number of the component $4$ between the components $2$ and $3$.

--3. An arbitrary well-defined homomorphism $Br_4 \to \Z$, which is a function of the windings numbers between components
is a linear combination of  $Lk(f)$ and $Lk(\tilde{f})$. 
\end{corollary}

Corollary $(\ref{winding})$ motivates the following definition.

\begin{definition}
Let $f \in Br_4$ be a (ordered) spherical braid. Define the total linking number $LK(f) \in \Z \oplus \Z$ by the following formula:
$$ LK(f) = Lk(f) \oplus Lk(\tilde f). $$
The total linking number is a well-defined homomorphism
$$ LK: Br_4 \mapsto \Z \oplus \Z. $$
\end{definition} 

\section{Hopf invariant of braids}

Let  $f \in Br_4$ be a (ordered) spherical braid with the trivial total linking number: $LK(f)=0$.
Such braids generate the subgroup in the group $Br_4$, denote this subgroup by $Brunn_4 \subset Br_4$.
Let us remark that this subgroup coincides no with the subgroup of Brunnian braids $Brun_4$, defined in \cite{B-M-V-W}, sequence 1.1.

\begin{theorem}\label{hopf}
There exists a well-defined homomorphism
\begin{eqnarray}\label{Hopf}
H: Brunn_4 \to \Z, 
\end{eqnarray}
called the Hopf invariant. The homomorphism $(\ref{Hopf})$ is invariant with respect to the action $(\ref{action})$ by an arbitrary permutation, which  in the kernel of $\theta$ in $(\ref{epi})$,
and is skew-invariant with respect to the action by a permutation, which is in the kernel of $\theta_1$ (the composition of $\theta$ with the projection on the first factor), but not in the kernel of $\theta_2$  (the composition of $\theta$ with the projection on the second factor).
\end{theorem}

\begin{remark}
Higher homotopy groups are described from the spherical  braids groups with non-ordered components in the \cite{B-M-V-W}, Sequence (1.1).
\end{remark}

\subsection*{Definition of the Hopf invariant}

Let $f \in Brun_4$ be an arbitrary. Consider the braid $f^{norm}$, given by $(\ref{Mobius})$.
Recal,  for the braid $f^{norm}$ the braid $g \in Br_3$, which consists of the straits $(1)$,$(2)$,$(3)$ of $f^{norm}$, is the constant braid at the points $0,1,\infty$ in $\hat \C$ correspondingly. 
Consider the strait (4) of the braid $f^{norm}$. This strait is represented by an oriented closed path  $i: S^1 \to \hat C \setminus \{0 \cup 1 \cup \infty\}$. This path determines  a cycle, which is an oriented boundary, because of the condition $LK(f^{norm})=0$. (Evidently, $LK(f^{norm})=LK(f)$, because the group of M$\ddot{\rm{o}}$bius transformations is connected.) 

Consider the inclusions
$$ 
I_0:  \hat\C \setminus \{0 \cup 1 \cup \infty \} \subset \hat\C \setminus \{1 \cup \infty\},
$$
$$ 
I_{\infty}:  \hat\C \setminus \{0 \cup 1 \cup \infty\} \subset \hat\C \setminus \{0 \cup 1\},
$$
$$ 
I_1: \hat\C \setminus \{0 \cup 1 \cup \infty\} \subset \hat\C \setminus \{0 \cup \infty\}.
$$
Because $H_1( \hat\C \setminus \{1 \cup \infty\};\Z) = \pi_1( \hat\C \setminus \{1 \cup \infty\})$,
for the homomorphism
$$I_{0,\sharp}: \pi_1( \hat\C \setminus \{0 \cup 1 \cup \infty\}) \to \pi_1( \hat\C \setminus \{0 \cup 1 \cup \infty\})$$
we get $I_{0,\sharp}([i])=0$. Analogously $I_{\infty,\sharp}([i])=0$, $I_{1,\sharp}([i])=0$.

There exist the following 3 maps of copies of the standard 2-disk
$$ e_{0}: D^2_0 \to \hat C \setminus \{1 \cup \infty\}, \quad e_0 \vert_{\partial D^2} = i,$$
$$ e_{\infty}: D^2_{\infty} \to \hat C \setminus \{0 \cup 1\}, \quad e_{\infty} \vert_{\partial D^2} = i,$$
$$ e_{1}: D^2_1 \to \hat C \setminus \{0 \cup \infty\}, \quad e_1 \vert_{\partial D^2} = i.$$

Consider a 2-sphere, which is represented by a gluing $D^2_0 \cup_{\partial} D^2_{\infty}$ of the disks
$D^2_0$, $D^2_{\infty}$ along the common boundary, which is identified with the circle $S^1_4$. Denote this sphere by $S^2_1$. Analogously define spheres 
$S^2_0=D^2_{\infty} \cup_{\partial} D^2_{1}$, $S^2_{\infty}=D^2_{1} \cup_{\partial} D^2_{0}$.

Consider the following commutative diagram of inclusions:
\begin{eqnarray}\label{0,1}
\begin{array}{ccc}
\hat \C \setminus \{0 \cup \infty \cup 1\} & \subset & \hat \C \setminus \{0 \cup \infty\} \\
\cap & & \cap \\
\hat \C \setminus \{\infty \cup 1\} & \subset & \hat \C \setminus \{\infty\} \\
\end{array}
\end{eqnarray}

Consider the mappings
$e_0: D_0^2 \to \hat \C \setminus \{1 \cup \infty\}$,
$e_1: D_1^2 \to \hat \C \setminus \{0 \cup \infty\}$
to the left bottom and to the right upper spaces of the diagram $(\ref{0,1})$ correspondingly.
The mapping $e_0 \cup_{\partial} e_1: S^2_{\infty} \to \hat \C \setminus \{\infty\}$
is well defined by gluing of the two mappings
$e_0$,
$e_1$
along the common mapping $i$ of the boundaries.  
Consider the standard 3-ball $D^3_{\infty}$ (with corners along the curve $S^1_4$) with the boundary $\partial D^3_{\infty} = S^2_{\infty}$. The mapping $e_0 \cup_{\partial} e_1$ can be extended to the mapping
\begin{eqnarray}\label{dinfty}
d_{\infty}: D^3_{\infty} \to \hat C \setminus \{\infty\}.
\end{eqnarray}
The target space of this mapping is the right bottom space of the diagram $(\ref{0,1})$. Because the target space of the mapping  $d_{\infty}$
is contractible, the mapping $d_{\infty}$ is well-defined up to homotopy. By the analogous constructions the following mappings
\begin{eqnarray}\label{d1}
d_{1}: D^3_{1} \to \hat \C \setminus \{1\},
\end{eqnarray}
\begin{eqnarray}\label{d0}
d_{0}: D^3_{0} \to \hat \C \setminus \{0\}
\end{eqnarray}
are well-defined.

The mappings $(\ref{dinfty})$, $(\ref{d1})$, $(\ref{d0})$ determine the mapping
\begin{eqnarray}\label{mapping}
h=h(f): S^3 \to S^2 
\end{eqnarray}
as follows. Take a 3-sphere $S^3$, which is catted into 3 balls $D^3_{\infty}, D^3_1, D^3_0$ along the common circle $S^1_4 \subset S^3$. 
The sphere $S^3$ is represented as the join $S^1_4 \ast S^1_a$ of the two standard circle. On the
circle $S^1_a$ take 3 points $x_0, x_{1}, x_{\infty} \in S^1_a$. The subsets $S^1_4 \ast [x_0,x_1] \subset S^3$,  
$S^1_4 \ast [x_1,x_{\infty}] \subset S^3$, $S^1_4 \ast [x_{\infty},x_0] \subset S^3$ are 3  copies of 3D disks, which are glued along corresponding subdomains in its boundaries. 

Let us identify $D^3_{\infty} \cong S^1_4 \ast [x_0,x_1]$, $D^3_{0} \cong S^1_4 \ast [x_1,x_{\infty}]$, $D^3_{1} \cong S^1_4 \ast [x_{\infty},x_0]$. The boundary $\partial D^3_{\infty}$ is identified with the balls $S^1_4 \ast \{0\} \cong D^2_0$, $S^1_4 \ast \{1\} \cong D^2_1$, which are glued along the common boundary $S^1_4$. The boundary $\partial D^3_{0}$ is identified with the balls $S^1_4 \ast \{1\} \cong D^2_1$, $S^1_4 \ast \{\infty\} \cong D^2_{\infty}$, which are identify along the same boundary $S^1_4$. The boundary $\partial D^3_{1}$ is identified with the balls $S^1_4 \ast \{\infty\} \cong D^2_{\infty}$, $S^1_4 \ast \{0\} \cong D^2_{0}$, which are identified along the same boundary $S^1_4$. 
The mappings $d_0$, $d_1$, $d_{\infty}$ on the corresponded balls are well-defined by the formulas
$(\ref{d0})$,$(\ref{d1})$,$(\ref{dinfty})$ correspondingly. This mappings define the mapping $(\ref{mapping})$ on the  3-sphere.

\begin{definition}\label{defHopf}
The Hopf invarian $H(f)$ for a braid $f \in Brunn_4$ in the formula $(\ref{Hopf})$ is defined as the Hopf invariant
of the mapping $h$ by the formula $(\ref{hopfinv})$. The mapping $h=h(f)$ is explicitly defined from the braid $f$ by the formula $(\ref{mapping})$.
\end{definition}

\subsection*{A formula to calculate  the Hopf invariant}

Let us introduces an explicit formula to calculate the Hopf invariant for a braid $f \in Brunn_4$.
Consider the complex plane $\C$. The 4-th strain of the braid $f^{norm}$ determines a curve on the plane without two points $\{0,1\}$, which will be denoted by
\begin{eqnarray}\label{pathgamma}
\gamma: S^1 \to \C \setminus \{0 \cup 1\}.
\end{eqnarray}
Let us consider the complex 1-form $(\ref{omega0})$. Define a complex 1-form 
\begin{eqnarray}\label{omega1}
\omega_1 =  \frac{1}{2 \pi \i} \frac{dz}{z-1}. 
\end{eqnarray}
Define a real (multivalued) function $\lambda_0$ by integration along the path $\gamma(t)$, $t \in [0,t] \subset S^1$ of the real part of the form $(\ref{omega0})$ as following:
\begin{eqnarray}\label{lambda0}
\lambda_0(t) = \Re {\int_0^t \omega_0}.
\end{eqnarray}
Define a real (multivalued) function $\lambda_1$ by integration along the path of the real part of the form $(\ref{omega1})$ as following:
\begin{eqnarray}\label{lambda1}
\lambda_1(t) = \Re {\int_0^t \omega_1}.
\end{eqnarray}
To take the multivalued functions $(\ref{lambda0})$, $(\ref{lambda1})$ well-defined, assume that the path $\gamma$ starts at the point $2 \in \C$: $\lambda_0(0)=2$,
$\lambda_1(0)=2$.

Define a closed 1-form $\psi(t)$ along a curve $\gamma(t) \in \C \setminus \{0 \cup 1\}$ by the following formula:
\begin{eqnarray}\label{form}
\psi(t) = \lambda_0(t) \omega_1 - \lambda_1(t) \omega_0.
\end{eqnarray}

Let us consider a function, which is well-defined as the real part of the integral
\begin{eqnarray}\label{int}
\Psi(T) = \Re{\int_{0}^T \psi(t) d\gamma}, \quad t \in [0,T] \subset S^1. 
\end{eqnarray}

\begin{theorem}\label{main}
The Hopf invariant of a braid $f \in Brun_4$ in the formula $(\ref{Hopf})$, which is defined by Definition $\ref{defHopf}$,
is calculated by the formula:
\begin{eqnarray}\label{m}
H(f)=\Psi(2\pi) = \frac{1}{2}\Re{ \int_{0}^{2\pi} \psi(t) d\gamma}, 
\end{eqnarray}
where $\gamma$ is the closed path, determined by the 4-th straight $S^1_4$ of the braid $f^{norm}$ by the formula $(\ref{pathgamma})$. 
\end{theorem}

From Theorem $\ref{main}$ we get a corollary.
\begin{corollary}\label{cor}

--1. The Hopf invariant $(\ref{Hopf})$ is an epimorphism.

--2. Assume that a braid $f \in Brunn_4$ is such that the braid $f^{norm}$ is represented by a commutator of the straight $(4)$ with straights $(1)$ and $(2)$ (such a braid is called the Borromean rings). 
Then $H(f)=\pm1$, where the sign in the formula depends on the sign of the commutator.
\end{corollary}

\subsubsection*{Proof of Corollary $\ref{cor}$}
It is sufficient to proof --2. The right side of the formula $(\ref{form})$ coincides with
the formula (28) \cite{B}, which is simplified for the considered example. For the Borromean ring the formula $\ref{m}$
is non-trivial. The factor $\frac{1}{2}$ in the right side of the formula provides $H(f)=1$ for the right
Borromean rings. Corollary is proved.

\section{Proof of Theorem $\ref{main}$}

Proofs of Lemma $\ref{linking}$, Corollary $\ref{winding}$, and of Theorem $\ref{hopf}$ are clear.
Let us proof Theorem $\ref{main}$. 
Consider the mapping $h: S^3 \to S^2 = \hat C$, which is defined by the formula $(\ref{mapping})$. Take volume forms $\Omega_0, \Omega_1 \in \Lambda^2(S^2)$, each form  is an ill--supported form at the point $0$, $1$ correspondingly.
The Hopf invariant $(\ref{Hopf})$ is calculated by the formula:
\begin{eqnarray}\label{integral}
H(f) = \frac{1}{2} \int_{S^3} h^{\ast}(\Omega_0) \wedge \beta_1 + h^{\ast}(\Omega_1) \wedge \beta_0, 
\end{eqnarray}
where $x \in S^3$, $h^{\ast}(\Omega_0) \in \Lambda^2(S^3)$ is the pull-back of $\Omega_0 \in \Lambda^2(S^2)$ by $h: S^3 \to S^2$,  $\beta_0 \in \Lambda^1(S^3)$ is an arbitrary 1-form, such that $d(\beta_0)=h^{\ast}(\Omega_0)$,
$\beta_1 \in \Lambda^1(S^3)$ is defined analogously to the definition of $\beta_0$. 

Evidently, the  1-forms $\beta_0$ in the integral $(\ref{integral})$ is represented up to a coboundary such that
$\beta_0 = 0$ inside the ball $D^3_0$. This follows from the fact that the curve $h^{-1}(0)$ is  outside the ball $D_0^3$.
Analogously, we may assume that  $\beta_1 = 0$ in the ball $D^3_1$.
Then we get the following simplification of $(\ref{integral})$:
$$H(f) = \frac{1}{2} \int_{D_{\infty}^3} h^{\ast}(\Omega_0) \wedge \beta_1 +  h^{\ast}(\Omega_1) \wedge \beta_0,$$
assuming that $\beta_1$ is inside $D_{\infty}^3 \cup D_0^3$, $\beta_0$ is inside $D_{\infty}^3 \cup D_1^3$.

In the ball $D_{\infty}^3 \cup D_0^3$ the 3-form $h^{\ast}(\Omega_0) \wedge \beta_1$ is exact, we get
$\alpha_1 \in \Lambda^2(D_{\infty}^3 \cup D_0^3)$, $d\alpha_1 = h^{\ast}(\Omega_0) \wedge \beta_1$.
Moreover, we may put $\alpha_1 = \beta_0 \wedge \beta_1$. 

In the ball $D_{\infty}^3 \cup D_1^3$ the 3-form $h^{\ast}(\Omega_1) \wedge \beta_0$ is exact, we get
$\alpha_0 \in \Lambda^2(D_{\infty}^3 \cup D_1^3)$, $d\alpha_0 = h^{\ast}(\Omega_1) \wedge \beta_0$.
Moreover, we may put $\alpha_0 = \beta_1 \wedge \beta_0 = -\beta_0 \wedge \beta_1$. 

Apply the 3D Gauss-Ostrogradsky formula, we get 
$$ \int_{D_{\infty}^3} h^{\ast}(\Omega_0) \wedge \beta_1  = \int_{D^2_1} \beta_0 \wedge \beta_1 , $$ 
$$ \int_{D_{\infty}^3}  h^{\ast}\Omega_1  \wedge \beta_0 = - \int_{D^2_{0}} \beta_0 \wedge \beta_1 . $$

The 2-form $\beta_0 \wedge \beta_1 \in \Lambda^2(D^2_1)$ is exact. Because in the disk $D^2_1$ the 0-form
$\lambda_1$ is well defined, and $d\lambda_1 = \beta_1$, we get:
$d(\lambda_1 \beta_0) = - \beta_0 \wedge \beta_1$.

Analogously, the 2-form $\beta_0 \wedge \beta_1 \in \Lambda^2(D^2_{0})$ is exact. Because in the disk $D^2_{0}$ the 0-form
$\lambda_0$ is well defined, and $d\lambda_0 = \beta_0$, we get:
$d(\lambda_0 \beta_1) = \beta_0 \wedge \beta_1$.

Apply the 2D Green  formula we get:
$$ \int_{D^2_1} \beta_0 \wedge \beta_1  = -\int_{\gamma} \lambda_1 \beta_0, $$
$$ \int_{D^2_0} \beta_0 \wedge \beta_1  = \int_{\gamma} \lambda_0 \beta_1, $$

The integral $(\ref{integral})$
is simplified as
$$ H(f) = \frac{1}{2} \int_{\gamma} \lambda_0 \beta_1 - \lambda_1 \beta_0. $$
This formula coincides with the formula $(\ref{m})$. Theorem $\ref{main}$ is proved.

\[  \]
\[  \]

Troitsk, Moscow region, IZMIRAN $\qquad \qquad \qquad \qquad$

 pmakhmet@mi.ras.ru  $\qquad \qquad \qquad \qquad$
\[  \] 

\end{document}